\documentclass[11pt,a4paper]{amsart}
\usepackage[T1]{fontenc}
\usepackage[utf8]{inputenc}
\usepackage{lmodern}
\usepackage[british]{babel}
\usepackage{amsmath,amsthm,amssymb,mathtools}
\usepackage[a4paper,margin=2.7cm]{geometry}
\usepackage{microtype}
\usepackage{enumitem}
\setlist[enumerate]{leftmargin=2.25em,itemsep=0.25em,topsep=0.35em}
\usepackage{xcolor}
\usepackage{tikz}
\usetikzlibrary{arrows.meta}
\providecommand{\figdir}{figures/}% Ordner der Abbildungen (TikZ-Quellen)
\definecolor{figPlus}{HTML}{4C4CFF}       % plus spins
\definecolor{figMinus}{HTML}{FF938D}      % minus spins, overwritten squares
\definecolor{figMinusLight}{HTML}{FFD3D0} % minus region beyond a square; squares not retained
\definecolor{figInterface}{HTML}{FF0000}  % interface line (schematic)
\definecolor{figPlusMark}{HTML}{0000FF}   % plus signs: spins frozen at +1
\definecolor{figBackground}{HTML}{F2FFF2} % background (green!5)
\definecolor{figDomain}{HTML}{BF0040}     % domain of the local chains; localisation region
\definecolor{figAccent}{HTML}{FF8000}     % paths
\definecolor{figGuide}{HTML}{333333}      % outlines, buffers, dimension marks

\definecolor{paperlink}{RGB}{0,0,114}
\usepackage{hyperref}
\hypersetup{colorlinks=true,linktoc=all,linkcolor=blue,citecolor=blue,urlcolor=blue,pdfstartview=Fit}
\usepackage{zref-clever}
\zcsetup{cap,nameinlink=true}

\theoremstyle{plain}
\newtheorem{theorem}{Theorem}[section]
\newtheorem{proposition}[theorem]{Proposition}
\newtheorem{lemma}[theorem]{Lemma}
\newtheorem{corollary}[theorem]{Corollary}
\theoremstyle{definition}

\theoremstyle{remark}

\AddToHook{env/proposition/begin}{\zcsetup{reftype=proposition}}
\AddToHook{env/lemma/begin}{\zcsetup{reftype=lemma}}
\AddToHook{env/corollary/begin}{\zcsetup{reftype=corollary}}
\AddToHook{env/definition/begin}{\zcsetup{reftype=definition}}
\AddToHook{env/assumption/begin}{\zcsetup{reftype=assumption}}
\AddToHook{env/example/begin}{\zcsetup{reftype=example}}
\AddToHook{env/remark/begin}{\zcsetup{reftype=remark}}
\zcRefTypeSetup{assumption}{Name-sg=Assumption,name-sg=assumption,Name-pl=Assumptions,name-pl=assumptions}
\numberwithin{equation}{section}

\IfFileExists{dsfont.sty}{%
  \usepackage{dsfont}%
}{%
  \DeclareMathAlphabet{\mathfplbb}{U}{fplmbb}{m}{n}%
}
\newcommand{\e}{\mathrm{e}}
\let\cedillaaccent\c
\DeclareRobustCommand{\c}{\ifmmode\mathrm{c}\else\expandafter\cedillaaccent\fi}
\let\dotunderaccent\d
\DeclareRobustCommand{\d}{\ifmmode\mathrm{d}\else\expandafter\dotunderaccent\fi}
\newcommand{\dd}{\,\mathrm{d}}
\newcommand{\TV}[1]{\lVert #1\rVert_{\mathrm{TV}}}
\DeclareMathOperator{\dist}{dist}

\hypersetup{pdftitle={Generic non-convergence of low-temperature Ising Glauber dynamics on Z2}}

\title[Generic non-convergence of Ising Glauber dynamics]{Generic non-convergence of low-temperature\\ Ising Glauber dynamics on~$\mathbb Z^2$}
\author[J. K\"oppl]{Jonas K\"oppl}
\author[Y. Steenbeck]{Yannic Steenbeck}

\address[Jonas K\"oppl]{Universität Münster,
Germany.}
\email{jonas.koeppl@uni-muenster.de}
\address[Yannic Steenbeck]{TU Braunschweig, Institut für Mathematische Stochastik,
Germany.}
\email{yannic.steenbeck@tu-braunschweig.de}
\date{\today}

\begin{document}

\begin{abstract}
We study how typical convergence is, in the topological sense of Baire category, for the heat-bath Glauber dynamics of the nearest-neighbour, zero-field Ising model on~$\mathbb Z^2$ started from translation-invariant initial laws. At sufficiently low positive temperature, we show that the initial laws whose evolution has both the plus and the minus phase as weak subsequential limits form a dense $G_\delta$ subset of the space of translation-invariant laws with the weak topology. Consequently, the initial laws for which the dynamics converges in distribution form a meagre set, which gives a partial answer to an open problem posed by Liggett. In particular, a large initial magnetisation does not force convergence to the plus phase, even for spatially ergodic initial laws. The proof combines equilibrium and mixing estimates of Gheissari and Sly with a periodic grid of plus spins of vanishing density and Baire's theorem, and it applies to all reversible monotone Glauber dynamics for the Ising model.
\end{abstract}

\maketitle

\section{Introduction}\label{sec:intro}

We consider the heat-bath Glauber dynamics of the Ising model on~$\mathbb Z^2$: at inverse temperature $\beta>0$, the spin at site $x\in\mathbb Z^2$ of a configuration $\sigma\in\{\pm1\}^{\mathbb Z^2}$ flips at rate
\begin{equation}\label{eq:rates}
c(x,\sigma)\coloneqq\frac{1}{1+\exp\big(2\beta\sum_{y\sim x}\sigma_x\sigma_y\big)}\,,
\end{equation}
where $y\sim x$ means that $x$ and~$y$ are nearest neighbours. Equivalently, at rate one the spin at~$x$ is resampled from its conditional law under the Ising Gibbs measures given the spins at the neighbours of~$x$. These rates are reversible with respect to the Ising Gibbs measures. For $\beta$ larger than the critical value~$\beta_{\mathrm c}$, the Ising model has two phases, the plus phase and the minus phase, and both are stationary for the dynamics. We ask whether the law of the dynamics started from a translation-invariant initial distribution converges as time tends to infinity. By compactness, we know that there must be convergent subsequences, so it would suffice to show that all subsequential limits agree. 

Convergence is known for some classes of initial laws. For instance, Gheissari and Sly~\cite{GS26} recently showed that, at low temperature, this dynamics started from a product law with a sufficiently strong plus bias converges to the plus phase. For other interacting particle systems, such as the voter model and the contact process, convergence theorems identify the limit for all translation-invariant initial laws, and the complete convergence theorem for the contact process even covers all initial laws; see~\cite[Chapters~V and~VI]{Lig05} and~\cite{Lig99}.
One could therefore be led to conjecture that the Glauber dynamics converges for every translation-invariant initial law, or at least for a typical one.

In this note we study how common convergent initial laws are in the topological sense of Baire category. This point of view goes back to the work of Oxtoby and Ulam~\cite{OU41} in the context of dynamical systems, who showed that a generic volume-preserving homeomorphism is ergodic, and it has since been taken in many contexts, for example for diffeomorphisms by Bonatti, Crovisier and Wilkinson~\cite{BCW09} and by Avila, Crovisier and Wilkinson~\cite{ACW16}, and more recently for translation-invariant point processes by Björklund~\cite{Bjo26}. Our main result shows that, at sufficiently low positive temperature, the conjecture fails and non-convergence is generic: in the space of translation-invariant initial laws with the weak topology, the laws whose evolution has both phases as subsequential limits form a dense $G_\delta$ set. 
Consequently, the initial laws for which the dynamics converges form a meagre set, which gives a partial answer to a problem posed by Liggett~\cite[Chapter~IV, Problem~7]{Lig05}. This in particular complements the results in~\cite{GS26}: even a spatially ergodic initial law with magnetisation arbitrarily close to one need not converge to the plus phase (\zcref{cor:ergodic}), so some assumption on the spatial structure of the initial law, such as the independence (or at least sufficiently fast decay of correlations) in~\cite{GS26}, is needed.

Without translation invariance, non-convergence is easier to obtain. At low temperature, Fontes, Isopi and Newman~\cite{FIN99} sketch a construction of an initial configuration, constant on concentric annuli of rapidly growing radii with alternating signs, whose evolution comes close to both phases at arbitrarily large times. This construction relies heavily on the spatial inhomogeneity of the initial configuration and cannot be made translation invariant.

In our translation-invariant setting, the proof combines equilibrium and mixing estimates of Gheissari and Sly with a grid construction: overwriting the spins on a sparse periodic grid of thin corridors by plus spins changes a translation-invariant law only slightly, but drives the dynamics close to the plus phase at a suitable later time. In this sense, our main result is essentially a corollary of the work of Gheissari and Sly, albeit a somewhat paradoxical one: the estimates behind their convergence results also show that convergence fails for a generic translation-invariant initial law. The argument applies to all reversible monotone Glauber dynamics for the Ising model, not only to the heat-bath dynamics.

\subsection*{Organisation of the manuscript}

In \zcref[S]{sec:setting}, we introduce the precise setting and notation before stating the main results. In \zcref[S]{sec:strategy}, we prove the main results, assuming an estimate for a single square of minus spins in the plus phase. We prove this estimate in \zcref[S]{sec:proofs}, based on the equilibrium and mixing estimates of Gheissari and Sly.

\section{Setting and main results}\label{sec:setting}

Configurations are elements of $\Omega\coloneqq\{\pm 1\}^{\mathbb Z^2}$, which carries the product topology and the coordinatewise partial order. The rates~\eqref{eq:rates} are reversible with respect to the Ising Gibbs measures with formal weight $\exp(\beta\sum_{\{x,y\}\colon x\sim y}\sigma_x\sigma_y)$, where the sum runs over all nearest-neighbour edges. Let $(S^\beta_t)_{t\ge0}$ be the associated semigroup, and let $\nu S^\beta_t$ be the law at time~$t$ of the process with initial law~$\nu$. We denote by $\mathcal P_\theta$ the set of translation-invariant probability measures on~$\Omega$, equipped with the weak topology. It is compact and metrisable, and $S^\beta_t$ maps it into itself, that is, it preserves translation invariance. We write $\Rightarrow$ for weak convergence, $\preceq$ for stochastic domination, $\mathcal G^\beta_\theta$ for the set of translation-invariant Ising Gibbs measures at inverse temperature $\beta$, and $\mu^+_\beta$, $\mu^-_\beta$ for the plus and the minus phase.

For a probability measure~$\nu$ on~$\Omega$, we call $\mu$ a \emph{weak subsequential limit} of $(\nu S^\beta_t)_{t\ge0}$ if $\nu S^\beta_{t_k}\Rightarrow\mu$ for some sequence $t_k\to\infty$. If $\nu\in\mathcal P_\theta$, then every weak subsequential limit of $(\nu S^\beta_t)_{t\ge0}$ lies in~$\mathcal G^\beta_\theta$. This attractor property goes back to Holley~\cite{Hol71}, see also~\cite{JK19,JK23} for extensions to more general, possibly irreversible, interacting particle systems. It also follows from~\cite[Corollary~2.2 and (2.4)]{KS26}, by which the specific relative entropy $h(\nu S^\beta_t\mid\mu^+_\beta)$ tends to zero, since $h(\cdot\mid\mu^+_\beta)$ is lower semicontinuous on~$\mathcal P_\theta$ and vanishes exactly on~$\mathcal G^\beta_\theta$ by the Gibbs variational principle~\cite[Chapter~15]{Geo11}. For $\beta>\beta_{\mathrm c}$, the set $\mathcal G^\beta_\theta$ is the segment $\{\lambda\mu^+_\beta+(1-\lambda)\mu^-_\beta\colon\lambda\in[0,1]\}$~\cite{Aiz80,Hig81}. The question is therefore whether $(\nu S^\beta_t)_{t \geq 0}$ converges to a single point of this line segment or continues to oscillate.

Convergence is known in several special cases. If the Gibbs measure is unique, for instance for $\beta<\beta_{\mathrm c}$, then $(\nu S^\beta_t)_{t\geq0}$ converges for every $\nu\in\mathcal P_\theta$. For $\beta>\beta_{\mathrm c}$, all the elements in~$\mathcal G^\beta_\theta$ are stationary with respect to the dynamics, so the semigroup $(S^\beta_t)_{t\geq0}$ cannot be ergodic. By attractiveness, $\nu S^\beta_t\Rightarrow\mu^+_\beta$ whenever $\mu^+_\beta\preceq\nu$, for instance for the all-plus configuration, and symmetrically $\nu S^\beta_t\Rightarrow\mu^-_\beta$ whenever $\nu\preceq\mu^-_\beta$. For more information on attractive spin systems, see~\cite[Chapter~III]{Lig05}. In the extreme case of zero temperature, a spin flips at rate one if it disagrees with the majority of its neighbours, at a fixed rate $\alpha\in(0,1]$ in case of a tie, and not at all otherwise. For $\alpha=1/2$, these are the limits of the rates~\eqref{eq:rates} as $\beta\to\infty$. For these zero-temperature dynamics, Fontes, Schonmann and Sidoravicius~\cite{FSS02} showed that, in every dimension at least two, the process started from a product law with density of plus spins sufficiently close to one fixates at~$+1$, that is, every spin is eventually equal to~$+1$ forever. For $\alpha=1/2$, Morris~\cite{Mor11} showed that the threshold density for this fixation tends to~$1/2$ as the dimension tends to infinity.

Recall that a subset of a topological space is a \emph{$G_\delta$ set} if it is a countable intersection of open sets, and that a set is \emph{meagre} if it is contained in a countable union of closed sets with empty interior. The complement of a dense $G_\delta$ set is meagre. Since $\mathcal P_\theta$ is compact and metrisable, Baire's theorem~\cite{Oxt80} shows that a dense $G_\delta$ subset of~$\mathcal P_\theta$ is not meagre itself. A property that holds on a dense $G_\delta$ set is therefore called \emph{generic}.

\begin{theorem}\label{thm:main}
There is $\beta_0<\infty$ such that, for every $\beta>\beta_0$, the set
\begin{equation*}
\mathcal N_\beta\coloneqq\big\{\nu\in\mathcal P_\theta\colon\mu^+_\beta\textup{ and }\mu^-_\beta\textup{ are weak subsequential limits of }(\nu S^\beta_t)_{t\ge0}\big\}
\end{equation*}
is a dense $G_\delta$ subset of~$\mathcal P_\theta$.
\end{theorem}

The value of~$\beta_0$ chosen in \zcref{sec:strategy} satisfies $\beta_0\ge\beta_{\mathrm c}$, so that $\mu^+_\beta\ne\mu^-_\beta$ for $\beta>\beta_0$. If $\nu S^\beta_t$ converges weakly as $t\to\infty$, then it has a single weak subsequential limit, so $\nu\notin\mathcal N_\beta$ because $\mu^+_\beta\ne\mu^-_\beta$. Since the complement of the dense $G_\delta$ set~$\mathcal N_\beta$ is meagre, \zcref{thm:main} has the following consequence.

\begin{corollary}\label{cor:meagre}
Let $\beta>\beta_0$. The set of $\nu\in\mathcal P_\theta$ for which $\nu S^\beta_t$ converges weakly as $t\to\infty$ is meagre in~$\mathcal P_\theta$.
\end{corollary}

For $\nu\in\mathcal N_\beta$, the weak subsequential limits of $(\nu S^\beta_t)_{t\ge0}$ in fact fill all of~$\mathcal G^\beta_\theta$, as follows from the attractor property, the continuity of $t\mapsto\nu S^\beta_t$ and the intermediate value theorem applied to the magnetisation.

As a last resort, one could hope that the failure to converge is due to initial laws that are not spatially ergodic, where $\nu\in\mathcal P_\theta$ is called \emph{spatially ergodic} if $\nu(A)\in\{0,1\}$ for every translation-invariant event~$A$. However, the spatially ergodic laws also form a dense $G_\delta$ subset of~$\mathcal P_\theta$~\cite[Comment~(14.13)]{Geo11}, and Baire's theorem~\cite{Oxt80} yields the following corollary.

\begin{corollary}\label{cor:ergodic}
Let $\beta>\beta_0$. The set of spatially ergodic $\nu\in\mathcal P_\theta$ for which $\mu^+_\beta$ and $\mu^-_\beta$ are weak subsequential limits of $(\nu S^\beta_t)_{t\ge0}$ is a dense $G_\delta$ subset of~$\mathcal P_\theta$.
\end{corollary}

In particular, since the laws $\nu\in\mathcal P_\theta$ with $\int\sigma_0\,\nu(\d\sigma)>1-\varepsilon$ form a non-empty open subset of~$\mathcal P_\theta$ for every $\varepsilon>0$, there are spatially ergodic initial laws with magnetisation arbitrarily close to one whose evolution does not converge.

The proof of \zcref[S]{thm:main} uses only general properties of the rates~\eqref{eq:rates}. It applies, with minor changes, to all reversible monotone Glauber dynamics for the Ising model, that is, to all single-spin-flip dynamics with positive, translation-invariant nearest-neighbour rates that are attractive and reversible with respect to the Ising Gibbs measures; the mixing estimate of~\cite[Theorem~4]{GS26} transfers to these dynamics by a comparison of Dirichlet forms. In particular, \zcref{thm:main} and \zcref{cor:meagre} hold for all such dynamics.

\section{Proof strategy and proof of the main result}\label{sec:strategy}

Throughout \zcref{sec:strategy,sec:proofs} we fix $\beta>\beta_0\coloneqq\max\{\beta_1,\beta_{\mathrm c}\}$, where $\beta_1$ is the constant of \zcref{lem:gs}, and we drop the index~$\beta$ from the notation: we write $S_t$, $\mu^+$, $\mu^-$ and~$\mathcal N$ for $S^\beta_t$, $\mu^+_\beta$, $\mu^-_\beta$ and~$\mathcal N_\beta$. Constants $c,C\in(0,\infty)$ depend only on~$\beta$ and may change from line to line. Signs $\sharp\in\{+,-\}$ are identified with $\pm1$ where convenient.

\subsection{Reduction via Baire's theorem}\label{sec:baire}
Let $d_{\mathcal P}$ be a metric for the weak topology on~$\mathcal P_\theta$. For $\sharp\in\{+,-\}$ and $n,N\in\mathbb N$ let
\begin{equation*}
\mathcal U^\sharp_{n,N}\coloneqq\bigcup_{t\in[N,\infty)}\big\{\nu\in\mathcal P_\theta\colon d_{\mathcal P}(\nu S_t,\mu^\sharp)<1/n\big\}\,.
\end{equation*}
By definition, $\mu^\sharp$ is a weak subsequential limit of $(\nu S_t)_{t\ge0}$ if and only if $\nu\in\mathcal U^\sharp_{n,N}$ for all $n,N\in\mathbb N$. Hence we can write
\begin{equation}\label{eq:intersection}
\mathcal N=\bigcap_{\sharp\in\{+,-\}}\,\bigcap_{n,N\in\mathbb N}\mathcal U^\sharp_{n,N}\,.
\end{equation}
For fixed~$t$, the map $\nu\mapsto\nu S_t$ is continuous on~$\mathcal P_\theta$ by the Feller property of the dynamics~\cite[Chapter~I]{Lig05}. Therefore every $\mathcal U^\sharp_{n,N}$ is open, and $\mathcal N$ is a $G_\delta$ set. By Baire's theorem~\cite{Oxt80}, a countable intersection of open dense subsets of the compact metrisable space~$\mathcal P_\theta$ is dense. In view of~\eqref{eq:intersection}, it therefore suffices to show that every $\mathcal U^\sharp_{n,N}$ is dense. This holds if, for $\sharp\in\{+,-\}$ and every $\rho\in\mathcal P_\theta$, there are $\rho_\ell\in\mathcal P_\theta$ and times $T_\ell\to\infty$ such that
\begin{equation}\label{eq:goal}
\rho_\ell\Rightarrow\rho\qquad\text{and}\qquad\rho_\ell S_{T_\ell}\Rightarrow\mu^\sharp\qquad\text{as }\ell\to\infty\,.
\end{equation}
In words, we have to approximate every translation-invariant law by translation-invariant laws whose evolution is close to~$\mu^\sharp$ at some arbitrarily late time.

\subsection{Monotone coupling}\label{sec:graphical}
We first construct a coupling so that all processes below live on one probability space that carries an independent collection of Poisson processes of rate one indexed by the sites of~$\mathbb Z^2$, independent marks that are uniformly distributed on~$[0,1]$ attached to their points, and, independently of these, the initial configurations. For a process $\omega=(\omega(t))_{t\ge0}$ with values in~$\Omega$, we write $\omega_x(t)$ for the spin at site~$x$ at time~$t$. At a point of the Poisson process at site~$x$ with mark~$u$, the spin at~$x$ is set to~$+1$ if $u\le p_\beta(s)$ and to~$-1$ otherwise, where
\begin{equation*}
p_\beta(s)\coloneqq\frac{\e^{\beta s}}{\e^{\beta s}+\e^{-\beta s}}
\end{equation*}
and $s$ is the sum of the spins at the neighbours of~$x$ just before the update. Thus the spin at~$x$ is resampled at rate one from its conditional law given the spins at the neighbours of~$x$, and it flips at rate~$c(x,\sigma)$. By Harris' argument~\cite{Har72}, this defines the dynamics for all initial configurations simultaneously, and almost surely the spin at a given site and time is determined by the initial configuration and finitely many Poisson points. We also use versions of the dynamics in which the spins in a given set are \emph{frozen}, i.e., never updated.

Each update is increasing in the configuration, since $p_\beta(s)$ is increasing in~$s$. Hence the dynamics is attractive in the following sense. Let $\omega=(\omega(t))_{t\ge0}$ and $\omega'=(\omega'(t))_{t\ge0}$ be two processes with $\omega(0)\le\omega'(0)$ such that every site frozen in~$\omega$ but not in~$\omega'$ carries~$-1$, and every site frozen in~$\omega'$ but not in~$\omega$ carries~$+1$. Then $\omega(t)\le\omega'(t)$ for all $t\ge0$. In particular, $\nu\preceq\nu'$ implies $\nu S_t\preceq\nu'S_t$ for all $t\ge0$.

\subsection{The grid construction}\label{sec:grid}
We construct the approximations required in~\eqref{eq:goal} by overwriting the spins on a periodic grid of corridors by plus spins, respectively by minus spins, with the grid placed at a uniformly random position.

A \textit{square of side~$n$} is a translate of $[-n/2,n/2]^2\cap\mathbb Z^2$ by a vector in~$\mathbb Z^2$, as in~\cite{GS26}. For odd~$n$, it has $n$ sites on each side. For $\Delta\subset\mathbb Z^2$ and $r\ge0$ let
\begin{equation*}
\Delta^{+r}\coloneqq\{x\in\mathbb Z^2\colon\dist_\infty(x,\Delta)\le r\}\,.
\end{equation*}
In particular, if $\Delta$ is a square of side~$n$ and $r$ is an integer, then $\Delta^{+r}$ is the concentric square of side $n+2r$. We write
\begin{equation*}
\partial\Delta\coloneqq\{x\in\mathbb Z^2\setminus\Delta\colon x\sim y\text{ for some }y\in\Delta\}
\end{equation*}
for the outer boundary of~$\Delta$ and $\lvert\Delta\rvert$ for the number of sites in~$\Delta$. For a probability measure~$\rho$ on~$\Omega$ and $\Delta\subset\mathbb Z^2$, we write $\rho_\Delta$ for the marginal of~$\rho$ on $\{-1,+1\}^\Delta$, and $\TV{\cdot}$ denotes the total variation distance, $\TV{\rho-\rho'}=\sup_A\lvert\rho(A)-\rho'(A)\rvert$, where the supremum runs over all measurable sets~$A$.

Let $\ell\in\mathbb N$, let $\Lambda_0\coloneqq\{-\ell,\dots,\ell\}^2$ be the square of side $2\ell+1$ centred at the origin, and let $L\coloneqq2\ell+1+8\lceil \ell^{3/4}\rceil$. For $v\in\mathbb Z^2$, let $\mathcal R^{(\ell)}_v\coloneqq\{\Lambda_0+v+Li\colon i\in\mathbb Z^2\}$. This is a family of squares of side $2\ell+1$, arranged with period~$L$ and separated by corridors of width $8\lceil \ell^{3/4}\rceil$; let $\Gamma^{(\ell)}_v$ be the set of \textit{corridor sites}, that is, of sites outside all squares in~$\mathcal R^{(\ell)}_v$. The corridor width is small compared with~$\ell$, so overwriting the corridors changes a translation-invariant law only slightly. It is large compared with the fluctuations of an Ising interface of length~$\ell$, which are of order~$\ell^{1/2}$, and this allows the plus corridors to drive the system to the plus phase at a suitable time, even if every square starts entirely minus.

For $\sharp\in\{+,-\}$, let $\Phi^{\sharp,(\ell)}_v\colon\Omega\to\Omega$ set all spins in~$\Gamma^{(\ell)}_v$ to~$\sharp$ and leave the other spins unchanged, i.e.,
\begin{equation}
\Phi^{\sharp,(\ell)}_v(\sigma)_x\coloneqq
\begin{cases}
\sharp &\text{if } x\in \Gamma^{(\ell)}_v,\\
\sigma_x &\text{otherwise.}
\end{cases}
\end{equation}
For $\rho\in\mathcal P_\theta$ and $\sharp\in\{+,-\}$ let
\begin{equation*}
\rho^\sharp_\ell\coloneqq L^{-2}\sum_{v\in\{0,\dots,L-1\}^2}\rho\circ\big(\Phi^{\sharp,(\ell)}_v\big)^{-1}
\end{equation*}
be the law of a configuration with law~$\rho$ whose corridors, for a uniformly random position of the grid, are overwritten by~$\sharp$. Since $\Gamma^{(\ell)}_{v+z}=\Gamma^{(\ell)}_v+z$ and $\Gamma^{(\ell)}_{v+Lz}=\Gamma^{(\ell)}_v$ for all $z\in\mathbb Z^2$, the translation invariance of~$\rho$ implies that also $\rho^\sharp_\ell\in\mathcal P_\theta$.

\begin{proposition}\label{prop:approx}
For $\ell\in\mathbb N$ let
\begin{equation}\label{eq:scales}
t_\ell\coloneqq\exp\big(C_1(\log(9\ell))^3\big)\qquad\text{and}\qquad T_\ell\coloneqq \ell\,t_\ell\,,
\end{equation}
where $C_1$ is the constant~$C$ of \zcref{lem:gs}. Then, for every $\rho\in\mathcal P_\theta$ and $\sharp\in\{+,-\}$,
\begin{equation*}
\rho^\sharp_\ell\Rightarrow\rho\qquad\text{and}\qquad\rho^\sharp_\ell S_{T_\ell}\Rightarrow\mu^\sharp\qquad\text{as }\ell\to\infty\,.
\end{equation*}
\end{proposition}

By \zcref{lem:gs}(c), $t_\ell$ bounds the mixing time of the dynamics in squares of side at most~$9\ell$ with plus boundary condition; this is used in the proof of \zcref{prop:square}. Moreover, $T_\ell\to\infty$ as $\ell\to\infty$, and $\log T_\ell\le C(1+\log \ell)^3$.

\begin{proof}[Proof of \zcref{thm:main}]
Let $\sharp\in\{+,-\}$ and $\rho\in\mathcal P_\theta$. By \zcref{prop:approx}, condition~\eqref{eq:goal} holds for $\rho_\ell\coloneqq\rho^\sharp_\ell$ and the times~$T_\ell$ of~\eqref{eq:scales}. By \zcref{sec:baire}, $\mathcal N$ is therefore a dense $G_\delta$ subset of~$\mathcal P_\theta$.
\end{proof}

\subsection{Approximating translation-invariant laws}\label{sec:approx}
We now prove \zcref{prop:approx}, given a localisation estimate for the process in which all squares of the grid are overwritten by minus spins. Let $\eta=(\eta(t))_{t\geq 0}$ be the process with initial law~$\mu^+$, so in particular $\eta$ is stationary. For $B\subset\mathbb Z^2$, let $\xi^B = (\xi^B(t))_{t\geq0}$ be the process started from
\begin{equation}\label{eq:xiB}
    \xi^B_x(0) \coloneqq
    \begin{cases}
        -1 &\text{if } x\in B,\\
        \eta_x(0) &\text{otherwise},
    \end{cases}
\end{equation}
i.e., we take the initial condition of $\eta$ but replace the spins in~$B$ by $-1$. For $v\in\mathbb Z^2$, we write $\xi^{(v,\ell)}\coloneqq\xi^{\mathbb Z^2\setminus\Gamma^{(\ell)}_v}$ for the process in which the spins on all the squares in the family $\mathcal{R}^{(\ell)}_v$ are replaced by $-1$. The next result, which we prove in \zcref{sec:square}, shows that on a fixed finite window, $\xi^{(v,\ell)}$ and the stationary process~$\eta$ agree at time~$T_\ell$ with probability tending to one.

\begin{proposition}\label{prop:family}
Let $T_\ell$ be as in~\eqref{eq:scales}. For every finite $\Delta\subset\mathbb Z^2$,
\begin{equation*}
\lim_{\ell\to\infty}\sup_{v\in\mathbb Z^2}\mathbb P\big(\xi^{(v,\ell)}(T_\ell)\ne\eta(T_\ell)\textup{ somewhere on }\Delta\big)=0\,.
\end{equation*}
\end{proposition}

\begin{figure}[t]
\centering
% fig_sandwich_body.tex -- sandwich of initial configurations (Section 3.3, Proposition 3.1).
% Needs: tikz, colors.tex, \figdir; data from make_figure_data.py.
% Exact: the three panels show the same grid position; pointwise order
% left <= middle <= right (asserted by the script).
% Illustration: sigma i.i.d. symmetric; plus phase drawn i.i.d. with minus
% density 0.05; square side 14 and corridor width 2 pixels (not to scale).
\begin{tikzpicture}[font=\small]
\input{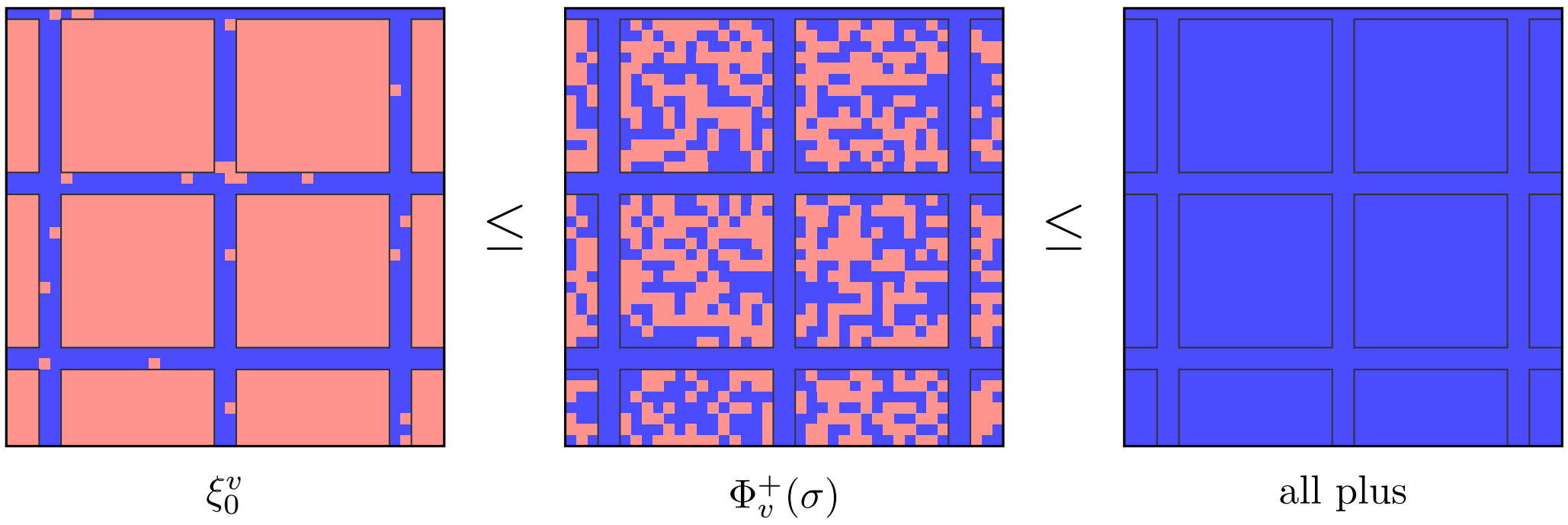}%
\def\px{0.095}%
\pgfmathsetmacro\W{\figSandN*\px}%
\def\gap{1.05}%
\def\figSandPanel#1#2{% #1: horizontal shift (cm), #2: data file (or empty)
  \begin{scope}[shift={({#1},0)}]
    \begin{scope}[x=\px cm, y=\px cm]
      \clip (0,0) rectangle (\figSandN,\figSandN);
      \fill[figPlus] (0,0) rectangle (\figSandN,\figSandN);
      \begin{scope}[fill=figMinus]
        \ifx&#2&\else\input{\figdir data/#2}\fi
      \end{scope}
      \foreach \k in {-1,...,2} {
        \foreach \m in {-1,...,2} {
          \draw[figGuide, thin] ({\figSandAX+\k*\figSandP},{\figSandAY+\m*\figSandP}) rectangle ++(\figSandL,\figSandL);
        }
      }
    \end{scope}
    \draw[semithick] (0,0) rectangle (\W,\W);
  \end{scope}}
\figSandPanel{0}{sandwich_left}
\figSandPanel{\W+\gap}{sandwich_mid}
\figSandPanel{2*\W+2*\gap}{}
\node[font=\Large] at ({\W+\gap/2},{\W/2}) {$\le$};
\node[font=\Large] at ({2*\W+3*\gap/2},{\W/2}) {$\le$};
\node[anchor=north] at ({\W/2},-0.12) {$\xi^v_0$};
\node[anchor=north] at ({3*\W/2+\gap},-0.12) {$\Phi^+_v(\sigma)$};
\node[anchor=north] at ({5*\W/2+2*\gap},-0.12) {all plus};
\end{tikzpicture}
\caption{The three initial configurations compared in the proof of \zcref{prop:approx}, for one position~$v$ of the grid; minus spins are red, plus spins blue, and the squares are outlined. Left: $\xi^{(v,\ell)}(0)$, which is~$-1$ on the squares and agrees with~$\eta(0)$, a configuration of the plus phase, on the corridors. Middle: $\Phi^{+,(\ell)}_v(\sigma)$, which agrees with~$\sigma$ on the squares and is~$+1$ on the corridors. Right: the all-plus configuration. The three configurations are ordered coordinatewise. They are illustrations, not to scale: $\sigma$ is i.i.d.\ symmetric, and the plus phase is drawn with a few isolated minus spins.}
\label{fig:sandwich}
\end{figure}

\begin{proof}[Proof of \zcref{prop:approx}]
By symmetry, it suffices to consider $\sharp=+$. Moreover, weak convergence of probability measures on~$\Omega$ is equivalent to the convergence of all finite-dimensional marginals, and these are determined by the probabilities of increasing events that depend on finitely many spins, by inclusion--exclusion over the events $\{\sigma_x=+1\text{ for all }x\in\Delta'\}$ with $\Delta'\subset\mathbb Z^2$ finite. For the second claim it therefore suffices to show that $\rho^+_\ell S_{T_\ell}(A)\to\mu^+(A)$ for every increasing event~$A$ that depends on finitely many spins.

\emph{Step 1: $\rho^+_\ell\Rightarrow\rho$.} Let $\Delta\subset\mathbb Z^2$ be finite. A configuration and its image under~$\Phi^{+,(\ell)}_v$ can differ on~$\Delta$ only if $\Delta$ intersects~$\Gamma^{(\ell)}_v$. For uniformly distributed $v\in\{0,\dots,L-1\}^2$, each site lies in~$\Gamma^{(\ell)}_v$ with probability
\begin{equation*}
1-\frac{(2\ell+1)^2}{L^2}\le\frac{2\big(L-(2\ell+1)\big)}{L}=\frac{16\lceil \ell^{3/4}\rceil}{L}\le\frac{8\lceil \ell^{3/4}\rceil}{\ell}\,.
\end{equation*}
Hence
\begin{equation*}
\TV{(\rho^+_\ell)_\Delta-\rho_\Delta}\le\frac{8\lvert\Delta\rvert\lceil \ell^{3/4}\rceil}{\ell}\longrightarrow0\,,
\end{equation*}
which proves $\rho^+_\ell\Rightarrow\rho$.

\emph{Step 2: $\rho^+_\ell S_{T_\ell}\Rightarrow\mu^+$.} Let $A$ be an increasing event that depends on finitely many spins. Fix~$v$, and let $\sigma$ have law~$\rho$, independently of everything else. Since $\xi^{(v,\ell)}(0)=-1$ on the squares of~$\mathcal R^{(\ell)}_v$ and $\Phi^{+,(\ell)}_v(\sigma)=+1$ on the corridors~$\Gamma^{(\ell)}_v$, we have $\xi^{(v,\ell)}(0)\le\Phi^{+,(\ell)}_v(\sigma)$, and both lie below the all-plus configuration; see \zcref{fig:sandwich}. Running the dynamics from these three initial configurations, attractiveness gives
\begin{equation*}
\mathbb P\big(\xi^{(v,\ell)}(T_\ell)\in A\big)\le\big(\rho\circ(\Phi^{+,(\ell)}_v)^{-1}\big)S_{T_\ell}(A)\le\delta_+S_{T_\ell}(A)\,,
\end{equation*}
where $\delta_+$ is the point mass at the all-plus configuration. Averaging over~$v$ gives
\begin{equation*}
L^{-2}\sum_{v\in\{0,\dots,L-1\}^2}\mathbb P\big(\xi^{(v,\ell)}(T_\ell)\in A\big)\le\rho^+_\ell S_{T_\ell}(A)\le\delta_+S_{T_\ell}(A)\,.
\end{equation*}
Since $\eta(T_\ell)$ has law~$\mu^+$, \zcref{prop:family} shows that the left-hand side converges to~$\mu^+(A)$. The right-hand side converges to~$\mu^+(A)$ as well, since $\delta_+S_t\Rightarrow\mu^+$ as $t\to\infty$, see~\cite[Chapters~III and~IV]{Lig05}. Hence $\rho^+_\ell S_{T_\ell}(A)\to\mu^+(A)$.
\end{proof}

\subsection{The single-square input}\label{sec:square}
The main input for \zcref{prop:family} is an estimate for a single square of minus spins in the plus phase. Let $\Lambda$ be a square of side $2\ell+1$, and set $\Lambda'\coloneqq\Lambda^{+\lceil \ell^{3/4}\rceil}$ and $\Lambda''\coloneqq\Lambda^{+3\lceil \ell^{3/4}\rceil}$, so that $\Lambda\subset\Lambda'\subset\Lambda''$. We consider the process~$\xi^\Lambda$ from~\eqref{eq:xiB}, in which only the spins in~$\Lambda$ are replaced by~$-1$.

\begin{proposition}\label{prop:square}
Let $T_\ell$ be as in~\eqref{eq:scales}. There are $c,C$ such that, for every $\ell\in\mathbb N$ and every square~$\Lambda$ of side $2\ell+1$, with probability at least $1-q_\ell$,
\begin{equation}\label{eq:square}
\xi^\Lambda(t)=\eta(t)\textup{ on }\mathbb Z^2\setminus\Lambda'\textup{ for all }t\in[0,T_\ell]\qquad\textup{and}\qquad\xi^\Lambda(T_\ell)=\eta(T_\ell)\,,
\end{equation}
where
\begin{equation*}
q_\ell\coloneqq C(1+T_\ell)\ell^C\e^{-c\sqrt \ell}+C\ell^22^{-\ell}\,.
\end{equation*}
\end{proposition}

In words, up to time~$T_\ell$ the processes $\xi^\Lambda$ and~$\eta$ differ only at sites of~$\Lambda'$, and at time~$T_\ell$ they agree everywhere. Since $\log T_\ell\le C(1+\log \ell)^3$, we have $T_\ell^kq_\ell\to0$ as $\ell\to\infty$ for every $k\in\mathbb N$.

\begin{figure}[tb]
\centering
% fig_family_body.tex -- a periodic family of squares (Section 3.4, Proposition 3.2).
% Needs: tikz with arrows.meta, colors.tex.
% Exact: corridor width 8 ceil(ell^{3/4}) and buffer width ceil(ell^{3/4}) in the ratio 8:1;
% period L = side + corridor width;
% buffers are disjoint; a square is drawn in the full minus colour (in \mathcal I)
% iff it meets Delta^{+kappa T_ell}. Schematic: ell^{3/4}/ell not to scale; kappa T_ell >> ell.
\begin{tikzpicture}[x=1cm,y=1cm,font=\small,
  dim/.style={figGuide,{Stealth[length=3.5pt]}-{Stealth[length=3.5pt]},thin}]
\def\s{1.5}% side of a square Lambda_i
\def\w{0.1}% buffer width ceil(ell^{3/4}); corridors are 8 times wider
\pgfmathsetmacro\p{\s+8*\w}
\def\Fx{6.5}\def\Fy{3.05}% centre of the window Delta (in a corridor)
\def\r{2.4}% half side of Delta^{+kappa T_ell} around Delta (schematic)
\fill[figBackground] (-0.4,-0.85) rectangle ({4*\p+\s+0.4},{2*\p+\s+0.65});
\foreach \i in {0,...,4} {
  \foreach \j in {0,...,2} {
    \pgfmathsetmacro\xa{\i*\p}
    \pgfmathsetmacro\ya{\j*\p}
    \pgfmathparse{and(and(\xa<=\Fx+\r,\xa+\s>=\Fx-\r),and(\ya<=\Fy+\r,\ya+\s>=\Fy-\r))}
    \ifnum\pgfmathresult=1
      \fill[figMinus] (\xa,\ya) rectangle ++(\s,\s);
    \else
      \fill[figMinusLight] (\xa,\ya) rectangle ++(\s,\s);
    \fi
    \draw[figGuide, densely dashed, very thin] ({\xa-\w},{\ya-\w}) rectangle ++({\s+2*\w},{\s+2*\w});
  }
}
% one square and its buffer, labelled
\node at ({\s/2},{2*\p+\s/2}) {$\Lambda_i$};
\node[figGuide, anchor=south west, inner sep=1.5pt] at ({-\w},{2*\p+\s+\w}) {$\Lambda'_i$};
% Delta^{+kappa T_ell}
\draw[figDomain, dashed, thick] ({\Fx-\r},{\Fy-\r}) rectangle ({\Fx+\r},{\Fy+\r});
% label in the horizontal corridor between the upper two rows, next to the right edge
\node[figDomain, anchor=west, inner sep=2pt] at ({\Fx+\r},{\p+\s+0.5*(\p-\s)}) {$\Delta^{+\kappa T_\ell}$};
% a nearest-neighbour path from a site x_0 outside Delta^{+kappa T_ell} to Delta
\draw[figAccent, thick, -{Stealth[length=4pt]}]
  (9.6,3.35) -- (8.75,3.35) -- (8.75,2.65) -- (7.6,2.65) -- (7.6,3.05) -- (6.66,3.05);
\fill[figAccent] (9.6,3.35) circle (1.6pt);
\node[anchor=north, inner sep=2pt] at (9.6,3.3) {$x_0$};
% the window Delta
\fill ({\Fx-0.07},{\Fy-0.07}) rectangle ({\Fx+0.07},{\Fy+0.07});
\node[anchor=south, inner sep=3pt] at (\Fx,{\Fy+0.05}) {$\Delta$};
% corridor width
\draw[dim] (\s,-0.25) -- node[below, font=\footnotesize] {$8\lceil\ell^{3/4}\rceil$} (\p,-0.25);
% period L of the grid
\draw[dim] ({2*\p},-0.25) -- node[below, font=\footnotesize] {$L$} ({3*\p},-0.25);
\end{tikzpicture}
\caption{The geometry in the proof of \zcref{prop:family} (schematic). Squares~$\Lambda_i$ of side $2\ell+1$, arranged with period~$L$, are separated by corridors of width $8\lceil \ell^{3/4}\rceil$, and their buffers $\Lambda'_i=\Lambda_i^{+\lceil \ell^{3/4}\rceil}$ (dashed) are disjoint. The squares that intersect~$\Delta^{+\kappa T_\ell}$ (dashed crimson) around the window~$\Delta$ are drawn darker, the other squares lighter; the indices of the darker squares form the set~$I$. The orange line is a nearest-neighbour path from a site~$x_0$ outside~$\Delta^{+\kappa T_\ell}$ to~$\Delta$. In the proof, $\kappa T_\ell$ is much larger than~$\ell$, so $I$ contains many more indices than drawn.}
\label{fig:family}
\end{figure}
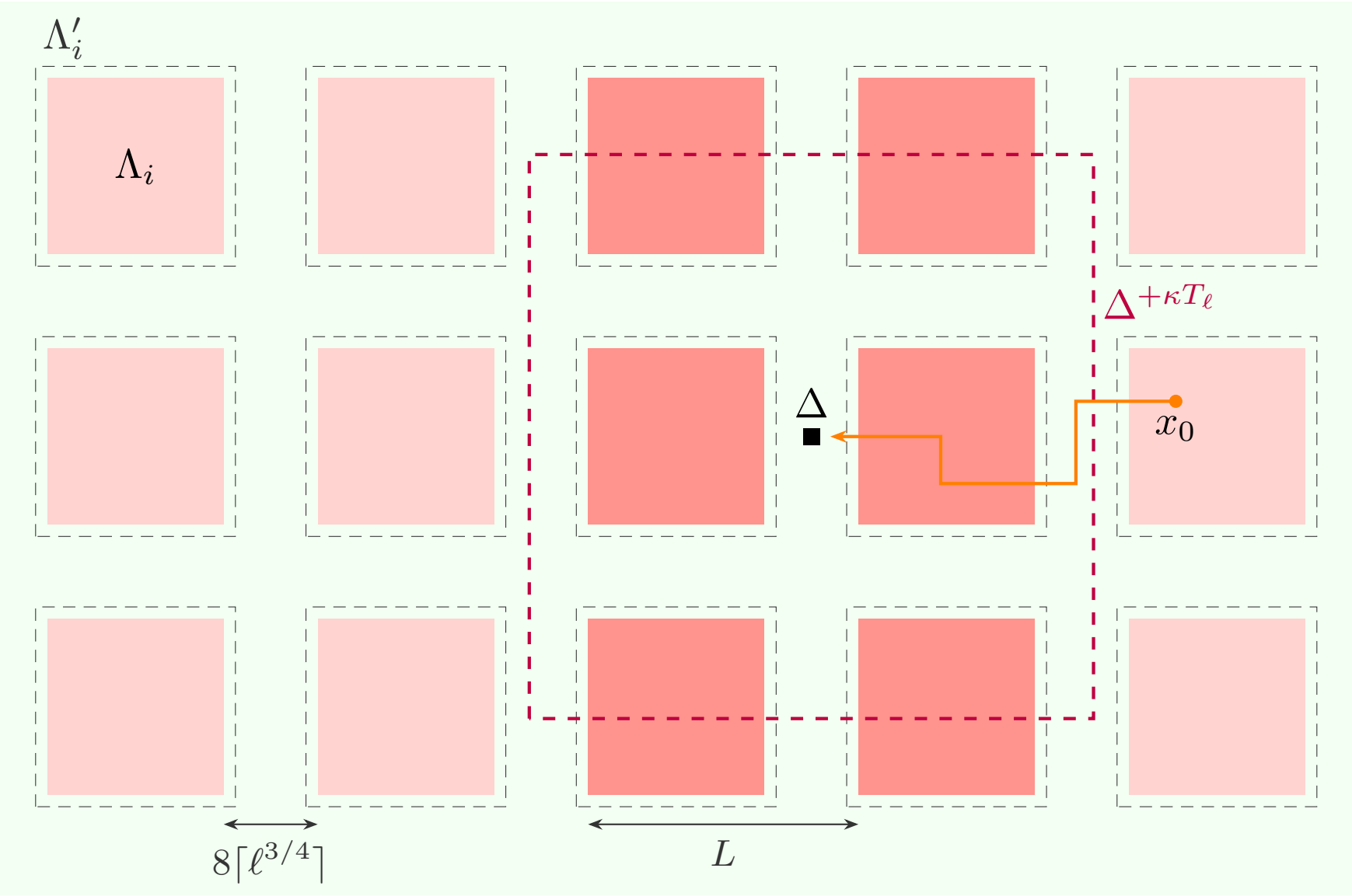

With this technical input at hand, we can complete the proof of \zcref{prop:family}; the argument is illustrated in \zcref{fig:family}.

\begin{proof}[Proof of \zcref{prop:family}]
Fix a finite set $\Delta\subset\mathbb Z^2$ and $v\in\mathbb Z^2$. For $i\in\mathbb Z^2$ let $\Lambda_i\coloneqq\Lambda_0+v+Li$, so that $\mathcal R^{(\ell)}_v=\{\Lambda_i\colon i\in\mathbb Z^2\}$, and let $\Lambda'_i\coloneqq\Lambda_i^{+\lceil \ell^{3/4}\rceil}$. Let $\kappa\coloneqq4\e^2$, let
\begin{equation*}
I\coloneqq\big\{i\in\mathbb Z^2\colon\Lambda_i\cap\Delta^{+\kappa T_\ell}\ne\emptyset\big\}
\end{equation*}
be the set of indices of the squares that intersect~$\Delta^{+\kappa T_\ell}$, and let $B_I\coloneqq\bigcup_{i\in I}\Lambda_i$. We first compare $\xi^{(v,\ell)}$ with~$\xi^{B_I}$, in which only the squares~$\Lambda_i$ with $i\in I$ are overwritten, and then $\xi^{B_I}$ with~$\eta$.

\emph{Step 1: Localisation.} The initial configurations of $\xi^{(v,\ell)}$ and~$\xi^{B_I}$ agree on~$\Delta^{+\kappa T_\ell}$. A disagreement between $\xi^{(v,\ell)}(T_\ell)$ and~$\xi^{B_I}(T_\ell)$ on~$\Delta$ therefore requires a nearest-neighbour path $x_0,\dots,x_j$ with $x_0\notin\Delta^{+\kappa T_\ell}$ and $x_j\in\Delta$, so that $j>\kappa T_\ell$, together with Poisson points at $x_1,\dots,x_j$ at increasing times in~$[0,T_\ell]$. For fixed~$j$, there are at most $\lvert\Delta\rvert4^j$ such paths, and for each of them the expected number of such $j$-tuples of Poisson points is $T_\ell^j/j!$. Since $j!\ge(j/\e)^j$ and $j>\kappa T_\ell=4\e^2T_\ell$, the expected number of such paths with Poisson points is at most
\begin{equation}\label{eq:paths}
\lvert\Delta\rvert\,4^j\,\frac{T_\ell^j}{j!}\le\lvert\Delta\rvert\Big(\frac{4\e T_\ell}{j}\Big)^j\le\lvert\Delta\rvert\,\e^{-j}\,.
\end{equation}
Summing~\eqref{eq:paths} over all integers $j>\kappa T_\ell$, we obtain
\begin{equation}\label{eq:localisation}
\mathbb P\big(\xi^{(v,\ell)}(T_\ell)\ne\xi^{B_I}(T_\ell)\text{ somewhere on }\Delta\big)\le\sum_{j>\kappa T_\ell}\lvert\Delta\rvert\,\e^{-j}=\frac{\lvert\Delta\rvert\,\e^{-\lfloor\kappa T_\ell\rfloor-1}}{1-\e^{-1}}\le2\lvert\Delta\rvert\,\e^{-\kappa T_\ell}\,.
\end{equation}

\emph{Step 2: Patching.} For $i\in I$, let $E_i$ be the event in~\eqref{eq:square} for $\Lambda=\Lambda_i$. By \zcref{prop:square}, $\mathbb P(E_i^{\c})\le q_\ell$ for every $i\in I$. On the event $\bigcap_{i\in I}E_i$, define $\tilde\xi(t)$ for $t\in[0,T_\ell]$ by $\tilde\xi_x(t)\coloneqq\xi^{\Lambda_i}_x(t)$ if $x\in\Lambda'_i$ for some $i\in I$, and $\tilde\xi_x(t)\coloneqq\eta_x(t)$ otherwise. Any two distinct sets~$\Lambda'_i$ are at distance at least $6\lceil \ell^{3/4}\rceil+1>2$ in the maximum norm, so the neighbourhood $\{x\}\cup\{y\colon y\sim x\}$ of a site~$x$ intersects at most one of them. If it intersects~$\Lambda'_i$ for some $i\in I$, then $\tilde\xi$ agrees with~$\xi^{\Lambda_i}$ on it: on~$\Lambda'_i$ by definition, and off~$\Lambda'_i$ because $\xi^{\Lambda_i}=\eta$ there on~$E_i$. If it intersects none of the sets~$\Lambda'_i$ with $i\in I$, then $\tilde\xi$ agrees with~$\eta$ on it. Hence $\tilde\xi$ follows the updates of the monotone coupling, and since $\tilde\xi(0)=\xi^{B_I}(0)$, it coincides with~$\xi^{B_I}$ on~$[0,T_\ell]$. On $\bigcap_{i\in I}E_i$ we therefore have $\xi^{B_I}(T_\ell)=\tilde\xi(T_\ell)=\eta(T_\ell)$, and a union bound gives
\begin{equation}\label{eq:finitefamily}
\mathbb P\big(\xi^{B_I}(T_\ell)\ne\eta(T_\ell)\big)\le\sum_{i\in I}\mathbb P(E_i^{\c})\le\lvert I\rvert\,q_\ell\,.
\end{equation}
This patching of single squares is related to~\cite[Observation~2.14 and Lemma~2.15]{GS26}.

\emph{Step 3: Counting.} The squares~$\Lambda_i$ with $i\in I$ are disjoint, contain $(2\ell+1)^2$ sites each and are contained in $\Delta^{+\kappa T_\ell+2\ell}$. Hence
\begin{equation}\label{eq:counting}
\lvert I\rvert\le\frac{\lvert\Delta^{+\kappa T_\ell+2\ell}\rvert}{(2\ell+1)^2}\le\lvert\Delta\rvert\,\frac{(2\kappa T_\ell+4\ell+1)^2}{(2\ell+1)^2}\le\lvert\Delta\rvert\Big(2+\frac{\kappa T_\ell}{\ell}\Big)^2\,.
\end{equation}

\emph{Step 4: Conclusion.} We have
\begin{align*}
&\mathbb P\big(\xi^{(v,\ell)}(T_\ell)\ne\eta(T_\ell)\text{ somewhere on }\Delta\big)\\
&\qquad\le\mathbb P\big(\xi^{(v,\ell)}(T_\ell)\ne\xi^{B_I}(T_\ell)\text{ somewhere on }\Delta\big)+\mathbb P\big(\xi^{B_I}(T_\ell)\ne\eta(T_\ell)\big)\\
&\qquad\le2\lvert\Delta\rvert\,\e^{-\kappa T_\ell}+\lvert I\rvert\,q_\ell\\
&\qquad\le2\lvert\Delta\rvert\,\e^{-\kappa T_\ell}+\lvert\Delta\rvert\Big(2+\frac{\kappa T_\ell}{\ell}\Big)^2q_\ell\,,
\end{align*}
where we used a union bound in the first step, \eqref{eq:localisation} and~\eqref{eq:finitefamily} in the second step, and \eqref{eq:counting} in the third step. The right-hand side does not depend on~$v$, and it tends to zero as $\ell\to\infty$, because $T_\ell\to\infty$ and $T_\ell^2q_\ell\to0$.
\end{proof}

\subsubsection*{Idea of the proof of \texorpdfstring{\zcref{prop:square}}{the proposition on one square}}
We compare $\xi^\Lambda$ and~$\eta$ with two stationary processes driven by the same updates: $\zeta$, frozen at~$-1$ on~$\Lambda$ and started from the plus phase conditioned to equal~$-1$ on~$\Lambda$, and~$\chi$, frozen at~$+1$ off~$\Lambda''$ and started from the Ising Gibbs measure in~$\Lambda''$ with plus boundary condition. Attractiveness gives $\zeta(t)\le\xi^\Lambda(t)\le\eta(t)\le\chi(t)$ for all $t\ge0$. The ring~$\partial\Lambda'$ lies at distance about~$\ell^{3/4}$ from~$\Lambda$ and at distance about~$2\ell^{3/4}$ from the complement of~$\Lambda''$. The equilibrium estimates of~\cite{GS26} therefore show that $\zeta$ and~$\chi$ agree on~$\partial\Lambda'$ at each fixed time with probability at least $1-C\ell^C\e^{-c\sqrt \ell}$, and counting the updates in~$\partial\Lambda'$ extends this to all times up to~$T_\ell$. Agreement on~$\partial\Lambda'$ forces $\xi^\Lambda=\eta$ off~$\Lambda'$. Inside~$\Lambda'$, both processes lie between the chains in~$\Lambda''$ started from all plus and from all minus, which coalesce by time~$T_\ell$ by the mixing estimate of~\cite{GS26}. \zcref[S]{fig:square} shows the geometry; the proof is given in \zcref{sec:proof-square}.

\begin{figure}[t]
\centering
% fig_square_body.tex -- one square (Section 3.4, Proposition 3.3).
% Needs: tikz with arrows.meta, colors.tex, \figdir (path prefix of figures/).
% Exact: Lambda in Lambda' = Lambda^{+ceil(ell^{3/4})} in Lambda'' = Lambda^{+3 ceil(ell^{3/4})};
% the ring is the outer boundary of Lambda'. Schematic: ell^{3/4}/ell not to scale; the
% interface (data/square_interface.tex) illustrates typical ell^{1/2} fluctuations, not a bound.
\begin{tikzpicture}[x=1cm,y=1cm,font=\small,
  dim/.style={figGuide,{Stealth[length=3.5pt]}-{Stealth[length=3.5pt]},thin},
  bar/.style={black,{Bar[width=4pt]}-{Bar[width=4pt]},semithick}]
\def\LH{2.2}% half side of Lambda (must match make_figure_data.py)
\def\w{0.75}% buffer width ceil(ell^{3/4})
\pgfmathsetmacro\H{\LH+\w}
\pgfmathsetmacro\Q{\LH+3*\w}
\pgfmathsetmacro\P{\Q+0.25}% position of the plus signs
% Lambda'': domain of the local chains chi and psi^\pm
\fill[figBackground] (-\Q,-\Q) rectangle (\Q,\Q);
\draw[figDomain, dashed, thick] (-\Q,-\Q) rectangle (\Q,\Q);
% spins frozen at +1 off Lambda''
\foreach \k in {0,...,24} {
  \pgfmathsetmacro\s{-\Q+0.1+\k*(2*\Q-0.2)/24}
  \node[figPlusMark, inner sep=0pt, font=\footnotesize] at (\s,\P) {$+$};
  \node[figPlusMark, inner sep=0pt, font=\footnotesize] at (\s,-\P) {$+$};
  \node[figPlusMark, inner sep=0pt, font=\footnotesize] at (\P,\s) {$+$};
  \node[figPlusMark, inner sep=0pt, font=\footnotesize] at (-\P,\s) {$+$};
}
% outer boundary of Lambda', where zeta and chi are compared
\draw[black, line width=1.8pt] ({-\H-0.06},{-\H-0.06}) rectangle ({\H+0.06},{\H+0.06});
% Lambda' = Lambda^{+ceil(ell^{3/4})}
\draw[figGuide, dashed] (-\H,-\H) rectangle (\H,\H);
% schematic interface around the minus square
\input{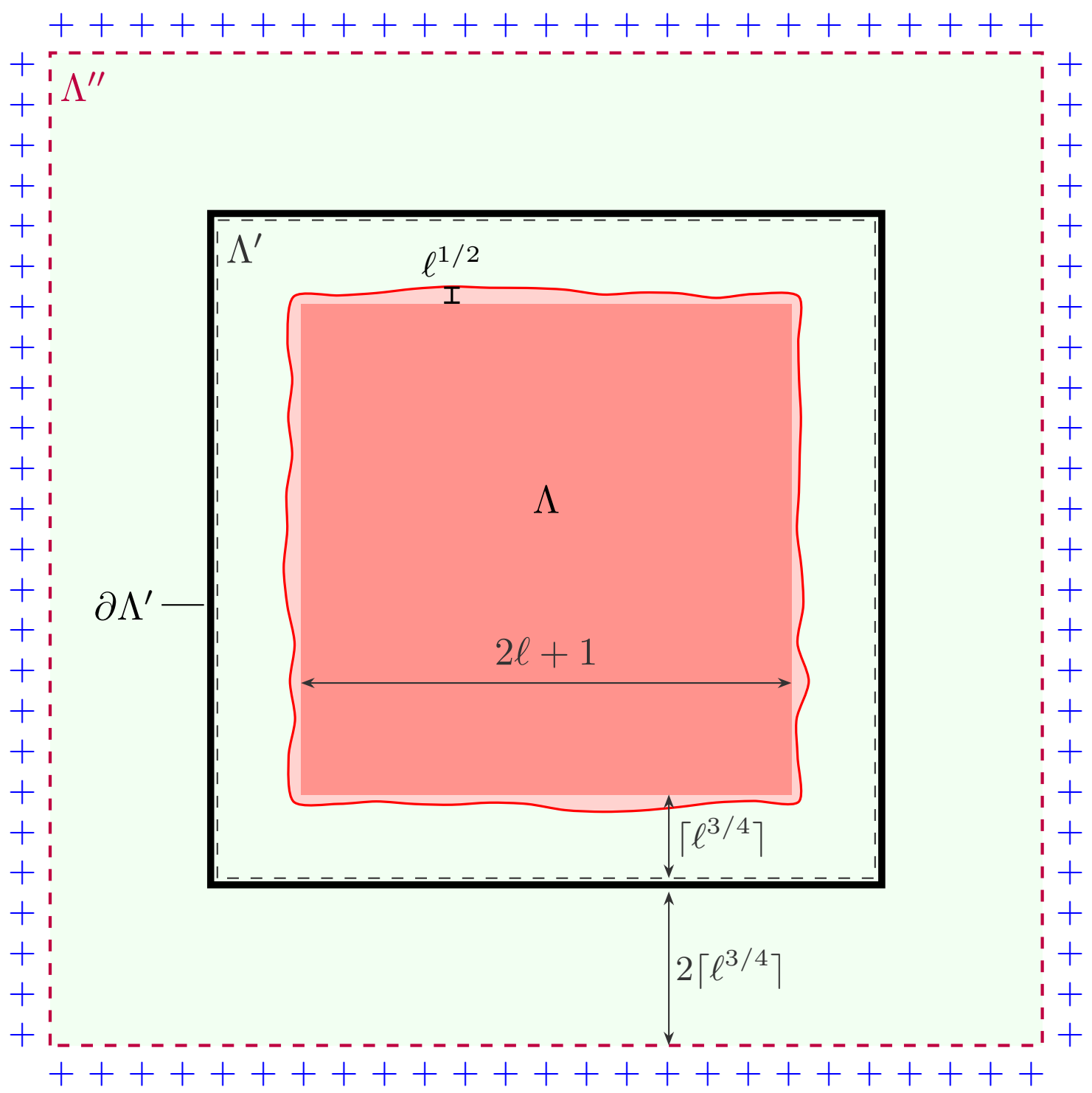}
% Lambda, overwritten by -1 (frozen at -1 in zeta)
\fill[figMinus] (-\LH,-\LH) rectangle (\LH,\LH);
% labels
\node at (0,0.45) {$\Lambda$};
\draw[dim] (-\LH,-1.2) -- node[above] {$2\ell+1$} (\LH,-1.2);
\node[figGuide] at ({-\H+0.25},{\H-0.25}) {$\Lambda'$};
\node[figDomain] at ({-\Q+0.3},{\Q-0.3}) {$\Lambda''$};
\draw[thin] ({-\H-0.12},-0.5) -- ({-\H-0.5},-0.5) node[left, inner sep=1.5pt] {$\partial\Lambda'$};
% vertical dimension markers below Lambda (the horizontal gaps are too narrow for the labels)
\draw[dim] (1.1,-\LH) -- node[right, font=\footnotesize, inner sep=1.5pt] {$\lceil\ell^{3/4}\rceil$} (1.1,-\H);
\draw[dim] (1.1,{-\H-0.12}) -- node[right, font=\footnotesize, inner sep=1.5pt] {$2\lceil\ell^{3/4}\rceil$} (1.1,-\Q);
\draw[bar] (\figSquareExcX,\LH) -- (\figSquareExcX,\figSquareExcY);
\node[font=\footnotesize, anchor=south, inner sep=1.5pt] at (\figSquareExcX,{\figSquareExcY+0.03}) {$\ell^{1/2}$};
\end{tikzpicture}
\caption{The geometry of \zcref{prop:square} (schematic, not to scale). The square~$\Lambda$ of side $2\ell+1$ is overwritten by~$-1$ (red); in the comparison process~$\zeta$ it is frozen at~$-1$, while $\chi$ and the extremal chains~$\psi^\pm$ live on $\Lambda''=\Lambda^{+3\lceil \ell^{3/4}\rceil}$ (dashed crimson) and are frozen at~$+1$ off~$\Lambda''$ (blue plus signs). The processes $\zeta$ and~$\chi$ are compared on the outer boundary~$\partial\Lambda'$ (thick black) of the buffer $\Lambda'=\Lambda^{+\lceil \ell^{3/4}\rceil}$ (dashed). In the equilibrium of~$\zeta$, an interface (drawn schematically) separates the minus spins attached to~$\Lambda$ from the plus phase; heuristically, its excursions are of order~$\ell^{1/2}$, much smaller than~$\ell^{3/4}$.}
\label{fig:square}
\end{figure}

\section{Proof of the single-square estimate}\label{sec:proofs}

\subsection{Gibbs measures and inputs from Gheissari and Sly}\label{sec:inputs}
The low-temperature input of this note consists of~\cite[Theorem~4 and Lemmas~3.3 and~3.5]{GS26}, which we recall in \zcref{lem:gs}; the rest of the argument is deduced from these estimates.

For finite $\Delta\subset\mathbb Z^2$ and $\tau\in\Omega$, let $\pi^\tau_\Delta$ be the Ising Gibbs measure in~$\Delta$ with boundary condition~$\tau$, viewed as a measure on~$\Omega$ concentrated on $\{\sigma\colon\sigma=\tau\text{ off }\Delta\}$. We write $\pi^+_\Delta$ for the plus boundary condition and, for $\Lambda\subset\Delta$, $\pi^\pm_{\Delta\setminus\Lambda}$ for the boundary condition that equals $-1$ on~$\Lambda$ and $+1$ off~$\Delta$. The measure $\pi^\tau_\Delta$ is stochastically increasing in~$\tau$~\cite[Chapter~3]{FV17}. Consequently, as $\Delta\uparrow\mathbb Z^2$, the measures $\pi^+_\Delta$ and $\pi^\pm_{\Delta\setminus\Lambda}$ converge weakly to~$\mu^+$ and to a measure~$\mu^\pm_\Lambda$ concentrated on $\{\sigma\colon\sigma=-1\text{ on }\Lambda\}$, respectively. Since conditioning~$\pi^+_\Delta$ on the event $\{\sigma=-1\text{ on }\Lambda\}$ yields~$\pi^\pm_{\Delta\setminus\Lambda}$, we have $\mu^\pm_\Lambda=\mu^+(\,\cdot\,\lvert\,\sigma=-1\text{ on }\Lambda)$, that is, $\mu^\pm_\Lambda$ is the plus phase conditioned to equal~$-1$ on~$\Lambda$. By the DLR equations and this monotonicity, $\mu^\pm_\Lambda\preceq\mu^+\preceq\pi^+_\Delta$ for all finite $\Lambda,\Delta\subset\mathbb Z^2$. By Strassen's theorem~\cite[Theorem~II.2.4]{Lig05}, these three measures admit a coupling that is almost surely ordered.

For a square~$\Delta$, let $\Omega^+_\Delta\coloneqq\{\sigma\in\Omega\colon\sigma_x=+1\text{ for all }x\in\mathbb Z^2\setminus\Delta\}$. The dynamics in~$\Delta$ frozen at~$+1$ off~$\Delta$ is the heat-bath dynamics in~$\Delta$ with plus boundary condition: it resamples the spin at each $x\in\Delta$ at rate one from its conditional law under~$\pi^+_\Delta$ given the other spins. It is a reversible Markov chain $\omega=(\omega(t))_{t\ge0}$ on the finite set~$\Omega^+_\Delta$ with stationary law~$\pi^+_\Delta$. Let
\begin{equation*}
d(t)\coloneqq\max_{\sigma\in\Omega^+_\Delta}\TV{\mathbb P_\sigma(\omega(t)\in\cdot\,)-\pi^+_\Delta}\qquad\text{and}\qquad t_{\mathrm{mix}}\coloneqq\inf\{t\ge0\colon d(t)\le1/4\}\,,
\end{equation*}
where $\mathbb P_\sigma$ refers to the chain started from~$\sigma$. Moreover, let
\begin{equation*}
\bar d(t)\coloneqq\max_{\sigma,\sigma'\in\Omega^+_\Delta}\TV{\mathbb P_\sigma(\omega(t)\in\cdot\,)-\mathbb P_{\sigma'}(\omega(t)\in\cdot\,)}\,.
\end{equation*}
Since $d\le\bar d\le2d$ and $\bar d$ is submultiplicative~\cite[Section~4.5]{LP17}, where the proofs apply verbatim in continuous time, we have
\begin{equation}\label{eq:submult}
d(kt)\le\bar d(kt)\le\bar d(t)^k\le\big(2d(t)\big)^k\le2^{-k}\qquad\text{for all }k\in\mathbb N\text{ and }t\ge t_{\mathrm{mix}}\,.
\end{equation}

\begin{lemma}\label{lem:gs}
There are $\beta_1,C<\infty$ such that the following holds for all $\beta>\beta_1$ and $n\ge1$.
\begin{enumerate}[label=(\alph*)]
\item \cite[Lemma~3.3]{GS26} Let $\Delta$ be a square of side~$n$, $0\le k\le n$, and $\Delta'$ the concentric square of side $n-k$. Then $\TV{(\pi^+_\Delta)_{\Delta'}-\mu^+_{\Delta'}}\le Cn\e^{-\beta k}$.
\item \cite[Lemma~3.5]{GS26} Let $\Lambda$ be a square of side~$n$, $0\le k\le n/100$, and $\Delta$ the concentric square of side $n+k$. Then $\TV{(\mu^\pm_\Lambda)_{\Delta'}-\mu^+_{\Delta'}}\le Cn^C\e^{-k^2/(Cn)}+C\e^{-\beta n}$ for every finite $\Delta'\subset\mathbb Z^2\setminus\Delta$.
\item \cite[Theorem~4]{GS26} The heat-bath dynamics in a square of side~$n$ with plus boundary condition satisfies $t_{\mathrm{mix}}\le\exp(C(1\vee\log n)^3)$.
\end{enumerate}
\end{lemma}

In~\cite{GS26}, part~(b) is stated for $\pi^\pm_{\Delta''\setminus\Lambda}$ and $\pi^+_{\Delta''}$, restricted to $\Delta''\setminus\Delta$, where $\Delta''\supset\Lambda$ is any concentric square of side at least~$n$; the bound is uniform in~$\Delta''$, and letting $\Delta''\uparrow\mathbb Z^2$ gives~(b).

\subsection{Proof of \texorpdfstring{\zcref{prop:square}}{the proposition on one square}}\label{sec:proof-square}
We may assume $2\lceil \ell^{3/4}\rceil\le(2\ell+1)/100$, since otherwise $q_\ell\ge1$ for a suitable choice of~$C$. Let $\zeta=(\zeta(t))_{t\ge0}$ be the process frozen at~$-1$ on~$\Lambda$ with initial law~$\mu^\pm_\Lambda$, and let $\chi=(\chi(t))_{t\ge0}$ be the process frozen at~$+1$ off~$\Lambda''$ with initial law~$\pi^+_{\Lambda''}$. By the ordered coupling of $\mu^\pm_\Lambda$, $\mu^+$ and~$\pi^+_{\Lambda''}$ from \zcref{sec:inputs}, we can choose the initial configurations such that $\zeta(0)\le\eta(0)\le\chi(0)$ holds almost surely. Since $\zeta(0)=-1$ on~$\Lambda$, also $\zeta(0)\le\xi^\Lambda(0)$, and attractiveness yields
\begin{equation}\label{eq:sandwich}
\zeta(t)\le\xi^\Lambda(t)\le\eta(t)\le\chi(t)\qquad\text{for all }t\ge0\,.
\end{equation}
Moreover, the process~$\chi$ is stationary because $\pi^+_{\Lambda''}$ is reversible for it. The process~$\zeta$ is stationary because $\mu^\pm_\Lambda$ is a Gibbs measure for the Ising specification on $\mathbb Z^2\setminus\Lambda$ with boundary condition~$-1$ on~$\Lambda$, and such Gibbs measures are reversible for the dynamics frozen on~$\Lambda$. This can easily be seen by checking that the proof for the unfrozen dynamics in~\cite[Chapter~IV]{Lig05} applies verbatim.

\emph{Step 1: $\zeta$ and~$\chi$ agree on~$\partial\Lambda'$ at each fixed time.} Let $x\in\partial\Lambda'$ and $t\ge0$. By~\eqref{eq:sandwich}, $\zeta_x(t)\le\chi_x(t)$. By the stationarity of $\zeta$ and of~$\chi$, the configurations $\zeta(t)$ and~$\chi(t)$ have laws $\mu^\pm_\Lambda$ and~$\pi^+_{\Lambda''}$, respectively. Only these individual laws enter here; we do not use the joint law of $\zeta(t)$ and~$\chi(t)$ beyond~\eqref{eq:sandwich}. Hence
\begin{equation}
\mathbb P(\zeta_x(t)\ne\chi_x(t))=\mathbb P(\chi_x(t)=+1)-\mathbb P(\zeta_x(t)=+1)=\pi^+_{\Lambda''}(\sigma_x=+1)-\mu^\pm_\Lambda(\sigma_x=+1)\,.
\end{equation}
We insert the plus phase and write the right-hand side as the sum of two differences,
\begin{equation}\label{eq:split}
\big(\pi^+_{\Lambda''}(\sigma_x=+1)-\mu^+(\sigma_x=+1)\big)+\big(\mu^+(\sigma_x=+1)-\mu^\pm_\Lambda(\sigma_x=+1)\big)\,,
\end{equation}
which we bound separately. The site~$x$ lies in $\Lambda^{+(\lceil \ell^{3/4}\rceil+1)}$, which is the square concentric with~$\Lambda''$ of side $(2\ell+1+6\lceil \ell^{3/4}\rceil)-(4\lceil \ell^{3/4}\rceil-2)$. Hence \zcref{lem:gs}(a) with $k=4\lceil \ell^{3/4}\rceil-2$ bounds the first difference in~\eqref{eq:split} by $C\ell\e^{-\beta \ell^{3/4}}$. The site~$x$ lies outside~$\Lambda'$, which is the square concentric with~$\Lambda$ of side $2\ell+1+2\lceil \ell^{3/4}\rceil$. Hence \zcref{lem:gs}(b) with $k=2\lceil \ell^{3/4}\rceil$ bounds the second difference in~\eqref{eq:split} by $C\ell^C\e^{-4\lceil \ell^{3/4}\rceil^2/(C(2\ell+1))}+C\e^{-\beta(2\ell+1)}$. Since $\ell^{3/4}\ge\sqrt \ell$ and $\lceil \ell^{3/4}\rceil^2/(2\ell+1)\ge\sqrt \ell/3$, both bounds are at most $C\ell^C\e^{-c\sqrt \ell}$. Since $\lvert\partial\Lambda'\rvert\le20\ell$, a union bound gives
\begin{equation}\label{eq:fixedtime}
\mathbb P(\zeta(t)\ne\chi(t)\text{ somewhere on }\partial\Lambda')\le C\ell^C\e^{-c\sqrt \ell}\qquad\text{for all }t\ge0\,.
\end{equation}

\emph{Step 2: $\zeta$ and~$\chi$ agree on~$\partial\Lambda'$ during a time interval.} Let $T>0$, and let $D(t)\in\{0,1\}$ indicate whether $\zeta(t)$ and~$\chi(t)$ differ somewhere on~$\partial\Lambda'$. If $D(t)=1$ for some $t\in[0,T]$, then $D(T)=1$ or $D$ jumps from~$1$ to~$0$ at some time in $(t,T]$. Such a jump happens at a Poisson point $(t',x)$ with $x\in\partial\Lambda'$, $t'\le T$ and $D(t'-)=1$. Since $D(t'-)$ is predictable, the expected number of Poisson points with these properties is $\lvert\partial\Lambda'\rvert\int_0^T\mathbb P(D(t')=1)\dd t'$. With~\eqref{eq:fixedtime} we obtain
\begin{equation}\label{eq:interval}
\begin{split}
\mathbb P(D(t)=1\text{ for some }t\in[0,T])&\le(1+\lvert\partial\Lambda'\rvert T)\sup_{t\in[0,T]}\mathbb P(D(t)=1)\\
&\le C(1+T)\ell^C\e^{-c\sqrt \ell}\,.
\end{split}
\end{equation}

\emph{Step 3: $\xi^\Lambda$ and~$\eta$ agree off~$\Lambda'$ up to time~$T_\ell$.} Let $G$ be the event that $\zeta(t)=\chi(t)$ on~$\partial\Lambda'$ for all $t\in[0,T_\ell]$. On~$G$, \eqref{eq:sandwich} gives $\xi^\Lambda(t)=\eta(t)$ on~$\partial\Lambda'$ for all $t\in[0,T_\ell]$. Moreover, $\xi^\Lambda(0)=\eta(0)$ off~$\Lambda$, and all neighbours of a site in $\mathbb Z^2\setminus(\Lambda'\cup\partial\Lambda')$ lie in $\mathbb Z^2\setminus\Lambda'$. By induction over the finitely many Poisson points that determine a given spin, at every update at such a site both processes see the same neighbouring spins, so the update is coupled perfectly. Hence $\xi^\Lambda(t)=\eta(t)$ on $\mathbb Z^2\setminus\Lambda'$ for all $t\in[0,T_\ell]$ on~$G$.

\emph{Step 4: $\xi^\Lambda$ and~$\eta$ agree on~$\Lambda'$ at time~$T_\ell$.} Let $\psi^+$ and~$\psi^-$ be the processes frozen at~$+1$ off~$\Lambda''$ and started from~$+1$ and from~$-1$ on~$\Lambda''$, respectively. Attractiveness gives $\psi^-(t)\le\chi(t)\le\psi^+(t)$ for all $t\ge0$, so on~$G$ both $\xi^\Lambda(t)$ and~$\eta(t)$ lie between $\psi^-(t)$ and~$\psi^+(t)$ on~$\partial\Lambda'$ for all $t\in[0,T_\ell]$. These inequalities hold on~$\Lambda'$ at time~$0$, and all neighbours of a site in~$\Lambda'$ lie in $\Lambda'\cup\partial\Lambda'$. By induction, at every update at such a site the processes see ordered neighbouring spins, and the update preserves this order because it is increasing in these spins. Hence the inequalities persist on~$\Lambda'$ up to time~$T_\ell$. Therefore $\xi^\Lambda(T_\ell)=\eta(T_\ell)$ on~$\Lambda'$ on the event $G\cap\{\psi^-(T_\ell)=\psi^+(T_\ell)\text{ on }\Lambda''\}$, and everywhere by Step~3.

\emph{Step 5: Probability of failure.} The square~$\Lambda''$ has side $n''\coloneqq2\ell+1+6\lceil \ell^{3/4}\rceil\le9\ell$, and $1\vee\log n''\le\log(9\ell)$ because $\log9>1$. Hence $t_\ell\ge t_{\mathrm{mix}}$ for the dynamics in~$\Lambda''$ by \zcref{lem:gs}(c), and $d(T_\ell)=d(\ell t_\ell)\le2^{-\ell}$ by~\eqref{eq:submult}. Since $\psi^-\le\psi^+$, for every $x\in\Lambda''$,
\begin{equation*}
\mathbb P(\psi^-_x(T_\ell)\ne\psi^+_x(T_\ell))=\mathbb P(\psi^+_x(T_\ell)=+1)-\mathbb P(\psi^-_x(T_\ell)=+1)\le2d(T_\ell)\le2^{1-\ell}\,.
\end{equation*}
A union bound over the at most $(9\ell)^2$ sites of~$\Lambda''$ gives
\begin{equation}\label{eq:coalescence}
\mathbb P\big(\psi^-(T_\ell)\ne\psi^+(T_\ell)\text{ somewhere on }\Lambda''\big)\le C\ell^22^{-\ell}\,.
\end{equation}
Let $F$ be the event in~\eqref{eq:square}. By Steps~3 and~4, $F\supset G\cap\{\psi^-(T_\ell)=\psi^+(T_\ell)\text{ on }\Lambda''\}$, and $G^{\c}$ is the event in~\eqref{eq:interval} for $T=T_\ell$. Hence
\begin{align*}
\mathbb P(F^{\c})&\le\mathbb P(G^{\c})+\mathbb P\big(\psi^-(T_\ell)\ne\psi^+(T_\ell)\text{ somewhere on }\Lambda''\big)\\
&\le C(1+T_\ell)\ell^C\e^{-c\sqrt \ell}+\mathbb P\big(\psi^-(T_\ell)\ne\psi^+(T_\ell)\text{ somewhere on }\Lambda''\big)\\
&\le C(1+T_\ell)\ell^C\e^{-c\sqrt \ell}+C\ell^22^{-\ell}=q_\ell\,,
\end{align*}
where we used Steps~3 and~4 in the first step, \eqref{eq:interval} with $T=T_\ell$ in the second step, and~\eqref{eq:coalescence} in the third step.\qed

\subsection*{Use of artificial intelligence}
The analysis of the proofs and results of \cite{GS26} and the consequences reported in this article were developed with substantial assistance of the AI systems ChatGPT Sol 6.1, Astra 6, and Astra 6.2 in an interactive process, while trying to understand the results of \cite{GS26} and what they imply for the question of convergence from non-product initial conditions and where proofs would break in that situation. The idea to use Baire's theorem to treat \textit{generic} initial conditions was inspired by the recent paper \cite{Bjo26}. An initial draft of the manuscript was then written by Claude Opus 5.5 under the direction of the authors. This initial version was then revised and rewritten to produce a legible version. During this process the authors carefully verified all of the mathematical details and take full responsibility for the correctness of the mathematical arguments in this note. 

\subsection*{Acknowledgements}
We are grateful for discussions about this question and related
problems on low-temperature Glauber dynamics with (in alphabetical order) Aernout van Enter, Nina Gantert, Reza Gheissari, Benedikt Jahnel, Christof Külske, Christian Maes, Fabio Martinelli, Charles Newman, Senya Shlosman, and Daniel Stein.

\subsection*{Funding}
The research of JK is funded by the Deutsche Forschungsgemeinschaft (DFG, German Research Foundation) under Germany's Excellence Strategy EXC 2044/2 –390685587, \textit{Mathematics Münster: Dynamics–Geometry–Structure}.

\bibliographystyle{amsalpha}
\bibliography{references}

\end{document}